\documentclass[a4paper,12pt]{amsart}
\usepackage{amsfonts}
\usepackage{amssymb}
\usepackage{ifthen}
\usepackage{graphicx}
\nonstopmode \numberwithin{equation}{section}
\usepackage[usenames]{color}
\newtheorem{thm}{Theorem}
\newtheorem{lem}{Lemma}
\newtheorem{cor}{Corollary}
\newtheorem{prop}{Proposition}
\newtheorem{ques}{Question}

\newtheorem{cl}{Claim}
\newtheorem{ca}{Case}
\newtheorem{sca}{Subcase}
\newtheorem{scl}{Subclaim}
\newtheorem{conj}[equation]{Conjecture}

\theoremstyle{definition}
\newtheorem{defn}{Definition}

\newtheorem{op}[equation]{Open Problem}
\newtheorem{rem}{Remark}
\newtheorem{exam}[equation]{Example}

\newcounter {own}
\def\theown {\thesection       .\arabic{own}}

\newenvironment{pf}[1][]{%
 \vskip 3mm
 \noindent
 \ifthenelse{\equal{#1}{}}%
  {{\slshape Proof. }}%
  {{\slshape #1.} }%
 }%
{\qed\bigskip}

\newcounter{alphabet}
\newcounter{tmp}
\newenvironment{Thm}[1][]{\refstepcounter{alphabet}%
\bigskip%
\noindent%
{\bf Theorem \Alph{alphabet}}%
\ifthenelse{\equal{#1}{}}{}{ (#1)}%
{\bf .} \itshape}{\vskip 8pt}

\makeatletter
\newcommand{\Ref}[1]{\@ifundefined{r@#1}{}{\setcounter{tmp}{\ref{#1}}\Alph{tmp}}}
\makeatother

\newcommand{\ID}{{\mathbb D}}

\def\be{\begin{equation}}
\def\ee{\end{equation}}

\newcommand{\bee}{\begin{enumerate}}
\newcommand{\eee}{\end{enumerate}}

\newcommand{\blem}{\begin{lem}}
\newcommand{\elem}{\end{lem}}
\newcommand{\bthm}{\begin{thm}}
\newcommand{\ethm}{\end{thm}}
\newcommand{\bcor}{\begin{cor}}
\newcommand{\ecor}{\end{cor}}
\newcommand{\beg}{\begin{exam}}
\newcommand{\eeg}{\end{exam}}
\newcommand{\begs}{\begin{examples}}
\newcommand{\eegs}{\end{examples}}
\newcommand{\bdefe}{\begin{defn}}
\newcommand{\edefe}{\end{defn}}
\newcommand{\bprob}{\begin{prob}}
\newcommand{\eprob}{\end{prob}}
\newcommand{\bques}{\begin{ques}}
\newcommand{\eques}{\end{ques}}
\newcommand{\bei}{\begin{itemize}}
\newcommand{\eei}{\end{itemize}}
\newcommand{\bcon}{\begin{conj}}
\newcommand{\econ}{\end{conj}}
\newcommand{\bop}{\begin{op}}
\newcommand{\eop}{\end{op}}

\newcommand{\bca}{\begin{ca}}
\newcommand{\eca}{\end{ca}}
\newcommand{\bsca}{\begin{sca}}
\newcommand{\esca}{\end{sca}}

\newcommand{\bcl}{\begin{cl}}
\newcommand{\ecl}{\end{cl}}

\newcommand{\bscl}{\begin{scl}}
\newcommand{\escl}{\end{scl}}

\newcommand{\bcons}{\begin{conjs}}
\newcommand{\econs}{\end{conjs}}
\newcommand{\bprop}{\begin{propo}}
\newcommand{\eprop}{\end{propo}}
\newcommand{\br}{\begin{rem}}
\newcommand{\er}{\end{rem}}
\newcommand{\brs}{\begin{rems}}
\newcommand{\ers}{\end{rems}}
\newcommand{\bo}{\begin{obser}}
\newcommand{\eo}{\end{obser}}
\newcommand{\bos}{\begin{obsers}}
\newcommand{\eos}{\end{obsers}}
\newcommand{\bpf}{\begin{pf}}
\newcommand{\epf}{\end{pf}}
\newcommand{\ba}{\begin{array}}
\newcommand{\ea}{\end{array}}
\newcommand{\beq}{\begin{eqnarray}}
\newcommand{\beqq}{\begin{eqnarray*}}
\newcommand{\eeq}{\end{eqnarray}}
\newcommand{\eeqq}{\end{eqnarray*}}

\newcounter{minutes}
\divide\time by 60
\newcounter{hours}
\multiply\time by 60 \addtocounter{minutes}{-\time}
\begin{document}
\bibliographystyle{amsplain}
\title [] {Geometric properties of $\log$-polyharmonic mappings}

\def\thefootnote{}
\footnotetext{ \texttt{\tiny File:~\jobname .tex,
          printed: \number\day-\number\month-\number\year,
          \thehours.\ifnum\theminutes<10{0}\fi\theminutes}
} \makeatletter\def\thefootnote{\@arabic\c@footnote}\makeatother

\author{Jiaolong Chen}
\address{Jiaolong Chen, Department of Mathematics,
Shantou University, Shantou, Guangdong 515063, People's Republic of China}
\email{jiaolongchen@sina.com}

\author{Bin Sheng}
\address{Bin Sheng, Department of Mathematics,
Shantou University, Shantou, Guangdong 515063, People's Republic of China}
\email{14bsheng@stu.edu.cn}

\author{Xiantao Wang${}^{~\mathbf{*}}$}
\address{Xiantao Wang, Department of Mathematics,
Shantou University, Shantou, Guangdong 515063, People's Republic of China}
\email{xtwang@stu.edu.cn}

\keywords{$\log$-polyharmonic, univalent, starlike, convex.\\
${}^{\mathbf{*}}$ Corresponding author}

\begin{abstract}
In this paper, a class of $\log$-polyharmonic mappings
$\mathcal{L}_p\mathcal{H}$ together with its subclass $\mathcal{L}_p\mathcal{H}(G)$ in the unit disk $\mathbb{D}=\{z: |z|<1\}$ is introduced, and
several geometrical properties such as the starlikeness, convexity
and univalence are investigated. In particular, we consider the Goodman-Saff  conjecture
and prove that the conjecture is true in $\mathcal{L}_p\mathcal{H}(G)$.
\end{abstract}

\thanks{
}

\maketitle \pagestyle{myheadings} \markboth{Jiaolong Chen, Bin Sheng and Xiantao Wang}{Geometric properties of $\log$-polyharmonic mappings}

\section{Introduction}\label{csw-sec1}
 A complex-valued mapping $f$ defined in a domain $D\subset
\mathbb{C}$ is called {\it harmonic} if $\Delta f=0$, where $\Delta$
represents the complex Laplacian operator
$$\Delta=4\frac{\partial^{2}}{\partial
z\partial\bar{z}}:=\frac{\partial^2}{\partial
x^2}+\frac{\partial^2}{\partial y^2}.$$ It is known that any
harmonic mapping $f$ in a simply connected domain $D$ can be written
in the form $f = g + \overline{h}$, where $g$ and $h$ are analytic
(cf. \cite{p}). A complex-valued mapping $F$ in a domain
$D$ is {\it biharmonic} if the Laplacian of $F$ is harmonic, that
is, $F$ satisfies the biharmonic equation $\Delta(\Delta F) = 0$. It
can be shown that in a simply connected domain $D$, every biharmonic
mapping has the form
$$F(z)=G_1(z)+|z|^2G_2(z),$$ where both $G_1$ and $G_2$ are harmonic in $D$.

More generally, a $2p$ $(p\geq 1)$ times continuously differentiable
complex-valued mapping $F$ in a domain $D \subset \mathbb{C}$ is {\it
polyharmonic} if $F$ satisfies the
polyharmonic equation $\Delta^{p}F =\Delta(\Delta^{p-1}F)= 0$. In a
simply connected domain $D$, a mapping $F$ is polyharmonic if and only
if $F$ has the following representation: \be\label{eq10} F(z)
=\sum_{k=1}^{p}|z|^{2(k-1)}G_{k}(z),\ee where each $G_{k}$ is
harmonic, i.e., $\Delta G_{k}(z) = 0$ for each $k\in \{1,\cdots,p\}$ (cf.
\cite{sh1, sh2}).

Obviously, when $p$ = 1 (resp. 2), $F$ is harmonic (resp. biharmonic).
The biharmonic equation arises when modeling and solving problems in a number of areas of science and engineering. Most importantly,
it is encountered in plane problems of elasticity and is used to describe creeping flows of viscous incompressible
 fluids. For instance, the determination of stress states around cavities in the stressed elastic body,
  regardless of cavity
 shapes, in its analytical approach has to be based on selection of a stress function that will satisfy the biharmonic equation
  (for details, see \cite{jh, s, w}).
However, the investigation of biharmonic mappings in the context of
geometric function theory is a recent one (cf. \cite{z1, z2, z3,
CW}). The reader is referred to \cite{ sh1, sh2, CRW1, CRW2, CRW3} for further discussions on polyharmonic mappings and \cite{sh5, jt, p}
as well as the references therein for the properties of harmonic mappings.

A mapping $G$ is said to be {\it log-harmonic} in a domain $D$ if
there is an analytic function $a$ such that $G$ is a solution of the
nonlinear elliptic partial differential equation
$$\frac{\overline{G_{\overline{z}}}}{\overline{G}}=a\frac{G_{z}}{G}.$$
It has been shown that if $G$ is a nonvanishing log-harmonic
mapping, then $G$ can be expressed as $G = f\overline{h}$, where
both $f$ and $h$ are analytic. A complex-valued mapping
$F$ in a domain $D$ is {\it log-biharmonic} if $\log F$ is
biharmonic, that is, the Laplacian of $\log F$ is harmonic. Further,
we say that $F$ is {\it $\log$-polyharmonic} if $\log F$ is
polyharmonic. It can be easily shown that every $\log$-polyharmonic
mapping $F$ in a simply connected domain $D$ has the form
\be\label{sun-1}F(z) = \prod_{k=1}^{p} G_{k}(z)^{|z|^{2(k-1)}},\ee
where all $G_{k}$ are nonvanishing log-harmonic mappings in $D$ for
$k\in \{1,\cdots,p\}$. When $p=1$ (resp. $p=2$), $\log$-polyharmonic
$F$ is $\log$-harmonic (resp. $\log$-biharmonic) (cf. \cite{lw1,
lw2}).

Some physical problems are modeled by $\log$-biharmonic mappings,
particularly, those arising in fluid flow theory and elasticity. The
log-biharmonic mappings are closely associated with the biharmonic
mappings, which appear in Stokes flow problems etc. Recently, the
properties of $\log$-harmonic and $\log$-biharmonic mappings have
been investigated by a number of authors (cf. \cite{z, zy, z4, sak, m}).

We say that a univalent polyharmonic mapping $F$ defined in $\mathbb{D}$, with $F (0) = 0$,
is {\it starlike} if the curve $F(re^{it})$ is starlike with respect
to the origin for each $0 < r < 1$. In other words, $F$ is starlike
if $$\frac{\partial}{\partial t}\big( \arg F(re^{it})\big)=\mbox{Re}
\frac{zF_{z}(z)-\overline{z}F_{\overline{z}}(z)}{F(z)} >0$$ for
$z\in \mathbb{D}\backslash\{0\}$.

A univalent  polyharmonic mapping $F$ defined in $\mathbb{D}$, with $\frac{\partial
}{\partial t}F(re^{it})\neq0$ whenever $0 < r < 1$, is said to be
{\it convex} if the curve $F (re^{it} )$ is convex for each $0 < r <
1$. In other words, $F$ is convex if $$ \frac{\partial} {\partial t}
\left(\arg   \frac{\partial}{\partial t}  F(re^{it}) \right)>0$$ for
$z\in \mathbb{D}\backslash\{0\}$.

Throughout this paper, we shall discuss $\log$-polyharmonic mappings defined
in $\mathbb{D}$.
We use $J_{F}$ to denote the Jacobian of $F$, that is,
$$J_{F} = |F_{z}|^{2}-|F_{\overline{z}}|^{2}.$$
It is known that $F$ is sense-preserving and locally univalent if $J_{F} > 0$. For convenience, we introduce the following notations:
 $$\mathcal{L}_p\mathcal{H}=\{F:\;F\;\mbox{has the expression \eqref{sun-1}}\;\mbox{in}\; \ID\}$$
 and
 $$\mathcal{L}_p\mathcal{H}(G) = \left\{F : F (z)=f(z)h(\overline{z})\prod _{k=1}^{p}G(z)^{\lambda_{k}|z|^{2(k-1)}}\right\},$$
 where $G$ is a non-vanishing $\log$-harmonic mapping in $\mathbb{D}$, $f$ and $h$ are non-vanishing analytic functions in $\mathbb{D}$,
 and all $\lambda_{k}$ are complex constants.

Obviously, we have the following.
\begin{prop}\label{prop1} Each element in $\mathcal{L}_p\mathcal{H}$ is $\log$-polyharmonic, and $\mathcal{L}_p\mathcal{H}(G)\subset \mathcal{L}_p\mathcal{H}$.\end{prop}

The main aim of this paper is two-fold. First, we discuss several geometrical properties such as the starlikeness, convexity and
 univalence for elements in $\mathcal{L}_p\mathcal{H}(G)$. Our results consist of four theorems which will be stated and proved in Section \ref{csw-sec4}. Second, we consider the Goodman-Saff conjecture (cf. \cite{sl}). Our result, Theorem \ref{thm5} in Section \ref{csw-sec5}, shows that the answer to this conjecture is positive in $\mathcal{L}_p\mathcal{H}(G)$. This result will be demonstrated in Section \ref{csw-sec5}.

\section{Properties of log-polyharmonic mappings}\label{csw-sec4}

First, we define a linear operator $\mathcal{L}$ by
$$\mathcal{L}=z\frac{\partial}{\partial z}-\overline{z}\frac{\partial}{\partial \overline{z}}.$$
The definition leads to the following two properties:

\begin{itemize}
  \item $\mathcal{L}[\alpha f + \beta g] = \alpha \mathcal{L}[f] + \beta \mathcal{L}[g]$,
  \item $\mathcal{L}[f g] = \mathcal{L}[f]g +f \mathcal{L}[g]$,
\end{itemize}
where $f$, $g$ are $C^{1}$ mappings and $\alpha$, $\beta$ are
complex constants.


\bthm\label{thm1} Suppose $F \in \mathcal{L}_p\mathcal{H}$ with
$F(z)=\prod_{k=1}^{p}G_{k}(z)^{|z|^{2(k-1)}}$. Then
\begin{enumerate}
  \item  $\mathcal{L}[\log F(z)]=\sum_{k=1}^{p}|z|^{2(k-1)}\mathcal{L}[\log G_{k}(z)]$;
  and
  \item $\mathcal{L}^{n}[\log F(z)]=\sum_{k=1}^{p}|z|^{2(k-1)}\mathcal{L}^{n}[\log G_{k}(z)]$,
\end{enumerate}
where $n\geq2$ is an integer.
\ethm

\bpf A simple calculation yields
$$\log F(z)=\sum_{k=1}^{p}|z|^{2(k-1)}\log G_{k}(z),$$
and for all $k\in \{1,\cdots,p\}$,
$$\mathcal{L}[|z|^{2(k-1)}]=0.$$
Using the product rule property and the linearity of the operator $\mathcal{L}$, we have
$$\mathcal{L}[\log F(z)]=\sum_{k=1}^{p}\mathcal{L}[|z|^{2(k-1)}\log G_{k}(z)]=\sum_{k=1}^{p}|z|^{2(k-1)}\mathcal{L}[\log G_{k}(z)],$$
which shows that the first statement in the theorem is true. The proof for the second statement easily follows from the first one
and mathematical induction.
\epf

By replacing each $G_k$ with $G^{\lambda_{k}}$ in Theorem
\ref{thm1}, we have the following corollary.

\bcor\label{cor1} Suppose $F\in \mathcal{L}_p\mathcal{H}(G)$ with $$F(z) =\prod_{k=1}^{p}
G(z)^{\lambda_{k}|z|^{2(k-1)}},$$ where
$\mathcal{L}[\log
G(z)]\neq0$ and $\sum_{k=1}^{p}\lambda_{k}|z|^{2(k-1)}\neq 0$. Then
$$ \frac{\mathcal{L}^{n}[\log F(z)]}{\mathcal{L}[\log F(z)]}= \frac{\mathcal{L}^{n}[\log G(z)]}{\mathcal{L}[\log G(z)]},\;\;n\geq2.$$
\ecor

From the definition of starlikeness, we know that $$F(z)
=\prod_{k=1}^{p} G(z)^{\lambda_{k}|z|^{2(k-1)}}\in \mathcal{L}_p\mathcal{H}(G)$$ is
starlike if and only if $F(0)=0$, $F$ is univalent and $\mbox{Re}
\frac{\mathcal{L}[F(z)]}{F(z)} >0$ for $z\in \mathbb{D}\backslash\{0\}$. Since $G$ is a non-vanishing
$\log$-harmonic mapping in $\mathbb{D}$, then it is impossible for
$F$ to reach this requirement. Instead, we consider the starlikeness
of $\log F$ and $\log G$.

\bthm\label{thm2}  Suppose that $F\in \mathcal{L}_p\mathcal{H}(G)$ with $$ F(z)
=\prod_{k=1}^{p} G(z)^{\lambda_{k}|z|^{2(k-1)}},$$ $\log F$ and $\log G$
are both univalent, and
$\sum_{k=1}^{p}\lambda_{k}|z|^{2(p-1)}\not=0$. Then
$\log G$ is starlike if and only if $\log F$ is starlike. \ethm \bpf
Obviously, for $z\in \mathbb{D}\backslash\{0\}$,
 \begin{align*}
\mbox{Re}\left( \frac{z(\log F(z))_{z}-\overline{z}(\log F(z))_{\overline{z}}}{\log F(z)}\right)
&=\mbox{Re}\frac{\sum_{k=1}^{p}\lambda_{k}|z|^{2(k-1)}\big(z(\log G(z) )_{z}- \overline{z } (\log G(z))_{\overline{z}}\big)}   {\sum_{k=1}^{p}\lambda_{k}|z|^{2(k-1)} \log G(z)}\\
&=\mbox{Re}\frac{z(\log G(z) )_{z}- \overline{z } (\log
G(z))_{\overline{z}}}   {\log G(z)}.
\end{align*} Then
the assumptions imply that $\log G$ is starlike if and only if $\log
F$ is starlike. \epf

To discuss the local univalence of $\log F$, where $F\in \mathcal{L}_p\mathcal{H}(G)$,
the following lemma is useful.

\blem\label{thm3} Suppose $F\in \mathcal{L}_p\mathcal{H}(G)$ with
$F(z)=f(z)h(\overline{z})\prod
_{k=1}^{p}G(z)^{\lambda_{k}|z|^{2(k-1)}}$. Then the Jacobian of
$\log F$, $J_{\log F}$, is given by
\begin{align*}
J_{\log F} (z)
=&\left| \frac{f'(z)}{f(z)} \right|^{2}-\left| \frac{h'(\overline{z})}{h(\overline{z})}
\right|^{2}+\left|B(z)\right|^{2}J_{\log G}(z)\\
 &+ 2|\log G(z)|^{2} {\rm Re} \left\{ \overline{A(z)}B(z) \mathcal{L}[\log (\log G(z))]\right\}\\
&+2 {\rm Re}\left\{ \overline{A(z)}\overline{\log G(z)} \mathcal{L}[\log
(f(z)h(\overline{z}))] \right\}+2{\rm Re} \left\{\overline{B(z)}
C(f,h, G; z)  \right \},
\end{align*}
where $$A(z)=\sum_{k=2}^{p}\lambda_{k}|z|^{2(k-2)}(k-1),\;
B(z)=\sum_{k=1}^{p}\lambda_{k}|z|^{2(k-1)}$$ and $$C(f,h, G;
z)=\frac{f'(z)}{f(z)}\frac{\overline{G_{z}(z)}}{\overline{G(z)}}-\frac{h'(\overline{z})}{h(\overline{z})}\frac{\overline{G_{\overline{z}}(z)}}{\overline{G(z)}}.$$
\elem
\bpf Taking the logarithm of $$F(z)
=f(z)h(\overline{z})\prod_{k=1}^{p} G(z)^{\lambda_{k}|z|^{2(k-1)}},$$
and then differentiating both sides with respect to $z$ and
$\overline{z}$, respectively, we get
$$\frac{F_{z}(z)}{F(z)}=\frac{f'(z)}{f(z)}+\overline{z}A(z)\log G(z)
+B(z)\frac{G_{z}(z)}{G(z)},$$
$$\frac{F_{\overline{z}}(z)}{F(z)}=\frac{h'(\overline{z})}{h(\overline{z})}+zA(z)\log G(z)
+B(z)\frac{G_{\overline{z}}(z)}{G(z)}.$$
Hence we have
\begin{align*}
\left|\frac{F_{z}(z)}{F(z)}\right|^{2}=&\left\{\frac{f'(z)}{f(z)}+\overline{z}A(z)\log G(z)
+B(z)\frac{G_{z}(z)}{G(z)} \right\}\cdot\left\{\frac{\overline{f'(z)}}{\overline{f(z)}}+z\overline{A(z)}\overline{\log G(z)}
+\overline{B(z)}\frac{\overline{G_{z}(z)}}{\overline{G(z)}} \right\}\\
=&\left|\frac{f'(z)}{f(z)} \right|^{2}+2\mbox{Re} \left\{  \frac{f'(z)}{f(z)} \left( z\overline{A(z)}\overline{\log G(z)}+\overline{B(z)}
\frac{\overline{G_{z}(z)}}{\overline{G(z)}} \right) \right\}\\
&+2{\rm
Re}\left\{z\overline{A(z)}B(z)
\overline{\log G(z)}\frac{G_{z}(z)}{G(z)}\right\}+|z|^{2}\left|A(z)\right|^{2}|\log G(z)|^{2}+\left|B(z)\right|^{2}\left|\frac{G_{z}(z)}{G(z)}\right|^{2}
\end{align*}
and
\begin{align*}
\left|\frac{F_{\overline{z}}(z)}{F(z)}\right|^{2}
=&\left|\frac{h'(\overline{z})}{h(\overline{z})}
\right|^{2}+2\mbox{Re}
\left\{\frac{h'(\overline{z})}{h(\overline{z})}
\left(\overline{z}\overline{A(z)}\overline{\log G(z)}+\overline{B(z)}\frac{\overline{G_{\overline{z}}(z)}}{\overline{G(z)}} \right) \right\}\\
&+2{\rm Re}\left\{\overline{z}\overline{A(z)}B(z)
\overline{\log G(z)}\frac{G_{\overline{z}}(z)}{G(z)}\right\}+|z|^{2}\left|A(z)\right|^{2}|\log G(z)|^{2}
+\left|B(z)\right|^{2}\left|\frac{G_{\overline{z}}(z)}{G(z)}\right|^{2}.
\end{align*}
Then, it follows from $ J_{\log F} (z)
=\frac{|F_{z}(z)|^{2}-|F_{\overline{z}}(z)|^{2}}{|F (z)|^{2}}$ that
\begin{align*}
&J_{\log F} (z)\\
=&\left| \frac{f'(z)}{f(z)} \right|^{2}-\left|\frac{h'(\overline{z})}{h(\overline{z})} \right|^{2}
+\left|B(z)\right|^{2}  \left(\left|\frac{G_{z}(z)}{G(z)}\right|^{2}-\left|\frac{G_{\overline{z}}(z)}{G(z)}\right|^{2}\right)\\
&+2{\rm Re}\left\{\overline{A(z)}B(z)\overline{\log G(z)}\left(    z\frac{G_{z}(z)}{G(z)}
-\overline{z}\frac{G_{\overline{z}}(z)}{G(z)}\right)\right\}\\
&+2 \mbox{Re}\left\{ \overline{A(z) }\left(\frac{zf'(z)}{f(z)}
-\frac{\overline{z}h'(\overline{z})}{h(\overline{z})}\right)\overline{\log G(z)}\right\}
+2\mbox{Re}\left\{\overline{B(z)}C(f,h,G;z) \right\}\\
=&\left| \frac{f'(z)}{f(z)} \right|^{2}-\left|\frac{h'(\overline{z})}{h(\overline{z})} \right|^{2}+\left|B(z)\right|^{2}
\left(\left|\frac{G_{z}(z)}{G(z)}\right|^{2}-\left|\frac{G_{\overline{z}}(z)}{G(z)}\right|^{2}\right)\\
&+2|\log G(z)|^{2}\mbox{Re}\left\{ \overline{A(z)}B(z) \frac{   zG_{z}(z)-\overline{z}G_{\overline{z}}(z)}{G(z)\log G(z)}\right\}\\
&+2\mbox{Re}\left\{  \overline{ A(z) }\left(\frac{zf'(z)}{f(z)}
-\frac{\overline{z}h'(\overline{z})}{h(\overline{z})}\right)\overline{\log G(z)}\right\}
+2\mbox{Re}\left\{\overline{B(z)}C(f,h,G;z) \right\}\\
=&\left| \frac{f'(z)}{f(z)} \right|^{2}-\left| \frac{h'(\overline{z})}{h(\overline{z})} \right|^{2}
+\left|B(z)\right|^{2}J_{\log G}(z)  + 2|\log G(z)|^{2} \mbox{Re} \left\{ \overline{A(z)}B(z) \mathcal{L}[\log (\log G(z))]\right\}\\
&+2 \mbox{Re}\left\{ \overline{A(z)}\overline{\log G(z)} \mathcal{L}[\log
(f(z)h(\overline{z}))] \right\} +2\mbox{Re}
\left\{\overline{B(z)}C(f,h,G;z) \right \}.
\end{align*} The proof is complete.
\epf

\bcor\label{cor9} Suppose $F\in \mathcal{L}_p\mathcal{H}(G)$ with
$F(z)=G(z)^{|z|^{2(p-1)}}$, where $p \geq 2$. Then the Jacobian of
$\log F$, $J_{\log F}$, is given by
\begin{align*}
J_{\log F} (z)=|z|^{4(p-1)}J_{\log G}(z)+ 2(p-1)|\log G(z)|^{2}
|z|^{2(2p-3)}{\rm Re} \left\{\mathcal{L}[\log (\log G(z))]\right\}.
\end{align*}
\ecor

\bthm\label{cor4} Suppose that $F\in \mathcal{L}_p\mathcal{H}(G)$ with
$$F(z)=f(z)h(\overline{z})\prod
_{k=1}^{p}G(z)^{\lambda_{k}|z|^{2(k-1)}},$$ the mapping $\log G$
is starlike and orientation preserving, $${\rm
Re}\Big(\overline{z}\frac{f'(\overline{z})}{f(\overline{z})}\mathcal{L}[\log
G(z)]\Big)>0,$$
$\frac{\overline{z}f'(\overline{z})}{f(\overline{z})}=\frac{zh'(z)}{h(z)}$,
and $\sum_{k=1}^{p}\lambda_{k}\neq0$ with $\lambda_{k}\geq0$
for all $k\in \{1,\cdots,p\}$. Then $\log F$ is orientation
preserving and consequently locally univalent for $z\in \mathbb{D}\backslash\{0\}$. \ethm

\bpf It follows from $$\frac{\overline{z}f'(\overline{z})}{f(\overline{z})}=\frac{zh'(z)}{h(z)}$$ that
$$\left| \frac{f'(z)}{f(z)} \right|^{2}=\left| \frac{h'(\overline{z})}{h(\overline{z})} \right|^{2}\;\; \mbox{and}\;\;
\mathcal{L}[\log (f(z)h(\overline{z}))]=\frac{zf'(z)}{f(z)}-\frac{\overline{z}h'(\overline{z})}{h(\overline{z})}=0.$$
Further, we get
\begin{align*}2\mbox{Re} \left( \frac{f'(z)}{f(z)}\frac{\overline{G_{z}(z)}}{\overline{G(z)}}
-\frac{h'(\overline{z})}{h(\overline{z})}\frac{\overline{G_{\overline{z}}(z)}}{\overline{G(z)}}
\right)&=\frac{2}{|z|^{2}}\mbox{Re} \left( \frac{zf'(z)}{f(z)}\frac{\overline{z}\overline{G_{z}(z)}}{\overline{G(z)}}
-\frac{\overline{z}h'(\overline{z})}{h(\overline{z})}\frac{z\overline{G_{\overline{z}}(z)}}{\overline{G(z)}} \right)\\
&=\frac{2}{|z|^{2}}\mbox{Re} \left(
\frac{\overline{z}f'(\overline{z})}{f(\overline{z})}\mathcal{L}[\log
G(z)]\right).
\end{align*}
Hence by Lemma \ref{thm3}, we see that
\begin{align*}
J_{\log F} (z) =&\left(\sum_{k=1}^{p}\lambda_{k}|z|^{2(k-1)}\right)\left\{\left(\sum_{k=1}^{p}\lambda_{k}
|z|^{2(k-1)}\right)J_{\log G}(z) +\frac{2}{|z|^{2}}\mbox{Re} \left( \frac{\overline{z}f'(\overline{z})}{f(\overline{z})}\mathcal{L}[\log G] \right) \right.\\
&  \left. + 2\left(\sum_{k=2}^{p}\lambda_{k}|z|^{2(k-2)}(k-1)\right)|\log G(z)|^{2} \mbox{Re} \big(\mathcal{L}[\log (\log G(z))]\big) \right \}.
\end{align*}

Since the assumptions imply that $\mbox{Re}(\mathcal{L}[\log(\log G)])
> 0$ for $z\in \mathbb{D}\backslash\{0\}$,
$J_{\log G} (z)
> 0$ and  $\sum_{k=1}^{p}\lambda_{k}|z|^{2(k-1)}>0$ for $z\in \mathbb{D}\backslash\{0\}$, it follows that $J_{\log F}  (z) > 0$ for $z\in \mathbb{D}\backslash\{0\}$, that
is, $\log F$ is orientation preserving, and hence $\log F$ is
locally univalent for $z\in \mathbb{D}\backslash\{0\}$. \epf

\bcor\label{cor10} Suppose that $F\in \mathcal{L}_p\mathcal{H}(G)$ with
$F(z)=G(z)^{|z|^{2(p-1)}}$, where $p\geq 2$, and that the mapping
$\log G$ is starlike and orientation preserving. Then $\log F$ is orientation preserving and consequently locally univalent
for $z\in \mathbb{D}\backslash\{0\}$. \ecor

We define the linear operator $\mathfrak{L}$ by
$$\mathfrak{L}=z\frac{\partial}{\partial z}+\overline{z}\frac{\partial}{\partial \overline{z}}.$$

\blem\label{thm4}  Suppose $F\in \mathcal{L}_p\mathcal{H}$ with $F(z)=\prod
_{k=1}^{p}G_{k}(z)^{|z|^{2(k-1)}}$. Then
\begin{enumerate}
 \item $-i\frac{\partial \log F(re^{it})}  {\partial t}=\sum_{k=1}^{p}|z|^{2(k-1)}L[\log G_{k}(z)],$
  \item $-\frac{\partial^{2} \log F(re^{it})}{\partial^{2} t}=\sum_{k=1}^{p}|z|^{2(k-1)}
  \big( \mathfrak{L}[\log G_{k}(z)]+z^{2}(\log G_{k}(z))_{zz}+\overline{z}^{2}(\log G_{k}(z))_{\overline{z}\overline{z}}\big).$
\end{enumerate}
\elem

\bpf Since
$$\log F(z)=\sum_{k=1}^{p}|z|^{2(k-1)}\log G_{k}(z),$$
it follows that
$$(\log F(z))_{z}=\sum_{k=2}^{p}(k-1)|z|^{2(k-2)}\overline{z}\log G_{k}(z)+\sum_{k=1}^{p}|z|^{2(k-1)}(\log G_{k}(z))_{z},$$
$$(\log F(z))_{\overline{z}}=\sum_{k=2}^{p}(k-1)|z|^{2(k-2)}z\log G_{k}(z)+\sum_{k=1}^{p}|z|^{2(k-1)}(\log G_{k}(z))_{\overline{z}}$$ and
\begin{align*}
(\log F(z))_{z\overline{z}}=&\sum_{k=2}^{p}(k-1)^{2}|z|^{2(k-2)}\log G_{k}(z)+\sum_{k=2}^{p}(k-1)|z|^{2(k-2)}\mathfrak{L}[\log G_{k}(z))].
\end{align*}
Therefore, by
\be\label{eq-4.2}
\frac{\partial \log F(re^{it})}{\partial t}
=iz(\log F(z))_{z}-i\overline{z}(\log F(z))_{\overline{z}},
\ee
we get
\begin{align*}
\frac{\partial \log F(re^{it})}{\partial t}
=&iz\left(\sum_{k=2}^{p}(k-1)|z|^{2(k-2)}\overline{z}\log G_{k}(z)+\sum_{k=1}^{p}|z|^{2(k-1)}(\log G_{k}(z))_{z}\right) \\
&-i\overline{z}\left(\sum_{k=2}^{p}(k-1)|z|^{2(k-2)}z\log G_{k}(z)+\sum_{k=1}^{p}|z|^{2(k-1)}(\log G_{k}(z))_{\overline{z}}\right) \\
=&i\sum_{k=1}^{p}|z|^{2(k-1)}L[\log G_{k}(z)],
\end{align*}from which the proof of part $(1)$ of the theorem follows.

Obviously, we know
\begin{align*}
&\mathfrak{L}[\log F(z)]-2|z|^{2}(\log F(z) )_{z\overline{z}}\\
=&2\sum_{k=3}^{p}(3k-k^{2}-2)|z|^{2(k-1)}\log G_{k}(z)+\sum_{k=1}^{p}(3-2k)|z|^{2(k-1)}   \mathfrak{L}[\log G_{k}(z)].\\
 \end{align*}
Upon differentiation we also have
\begin{align*}
(\log F(z) )_{zz} = &\sum_{k=3}^{p}(k-1)(k-2)|z|^{2(k-3)}\overline{z}^{2}
\log G_{k}(z)+2\sum_{k=2}^{p}(k-1)|z|^{2(k-2)}\overline{z}(\log G_{k}(z))_{z}\\
&+\sum_{k=1}^{p}|z|^{2(k-1)}(\log
G_{k}(z))_{zz}\end{align*} and
\begin{align*}(\log F(z)
)_{\overline{z}\overline{z }}=
&\sum_{k=3}^{p}(k-1)(k-2)|z|^{2(k-3)}z^{2}
\log G_{k}(z)+2\sum_{k=2}^{p}(k-1)|z|^{2(k-2)}z(\log G_{k}(z))_{\overline{z}}\\
&+\sum_{k=1}^{p}|z|^{2(k-1)}(\log
G_{k}(z))_{\overline{z}\overline{z}}.
 \end{align*}
Hence, we get
 \begin{align*}
&z^{2}(\log F (z))_{zz} + \overline{z}^{2}(\log F (z) )_{\overline{z}\overline{z }} \\
=  &2\sum_{k=3}^{p}(k-1)(k-2)|z|^{2(k-1)}\log G_{k} (z)+2\sum_{k=2}^{p}(k-1)|z|^{2(k-1)}\mathfrak{L}[\log G_{k}(z)] \\
&+\sum_{k=1}^{p}|z|^{2(k-1)}\big(z^{2}(\log G_{k} (z))_{zz}+\overline{z}^{2}(\log G_{k} (z))_{\overline{z}\overline{z}}\big).
 \end{align*}
We infer that
\beq\label{eq-4.3}
&\;\;&-\frac{\partial^{2} \log F(re^{it})}{\partial^{2} t}\\ \nonumber
&=&-\frac{\partial}{\partial t}[iz(\log F(z))_{z}-i\overline{z}(\log F(z))_{\overline{z}}]\\ \nonumber
&=& \mathfrak{L}[\log F(z)]-2|z|^{2}(\log F(z))_{z\overline{z}} +z^{2}(\log F(z))_{zz} +\overline{z}^{2}(\log F(z))_{\overline{zz}}\\ \nonumber
&=&\sum_{k=1}^{p}|z|^{2(k-1)}\big( \mathfrak{L}[\log G_{k}(z)]+z^{2}
  (\log G_{k}(z))_{zz}+\overline{z}^{2}(\log
  G_{k}(z))_{\overline{z}\overline{z}}\big),
\eeq
 from which the proof of part $(2)$ follows.
\epf

\bthm\label{cor5} Suppose that $F\in \mathcal{L}_p\mathcal{H}(G)$ with
$$F(z)=f(z)h(z)\prod _{k=1}^{p}G(z)^{\lambda_{k}|z|^{2(k-1)}},$$ and
suppose further that both $f$ and $h$ are constant functions, $\log G$
and $\log F$ are univalent, $\mathcal{L}[\log G(z)]\neq0$ and
$\sum_{k=1}^{p}\lambda_{k}|z|^{2(k-1)}\neq0$ for $z\in \mathbb{D}\backslash\{0\}$. Then
$\log G$ is convex if and only if $\log F$ is convex.\ethm
\bpf It follows from \eqref{eq-4.2} and \eqref{eq-4.3} in the proof of Lemma \ref{thm4} that for $z\in \mathbb{D}\backslash\{0\}$,
\begin{align*}
&\frac{\partial}{\partial t}\left(\arg \frac{\partial \log
F(re^{it})}{\partial t}\right)
=\mbox{Im} \left(\frac{\frac{\partial^{2}\log F(re^{it})} {\partial t^{2}}}{\frac{\partial\log F(re^{it})}{\partial t} } \right)\\
=&\mbox{Re} \left( \frac{  \mathfrak{L}[\log F(z)]-2|z|^{2}(\log F(z) )_{z\overline{z}} + z^{2}(\log F(z) )_{zz}
+ \overline{z}^{2}(\log F (z))_{\overline{z}\overline{z}}}  {\mathcal{L}[\log F(z)]} \right)\\
=&\mbox{Re} \left( \frac{\mathfrak{L}[\log G(z)]  -2|z|^{2}(\log G(z) )_{z\overline{z}} + z^{2}(\log G(z) )_{zz}
+ \overline{z}^{2}(\log G(z) )_{\overline{z}\overline{z}}}  {\mathcal{L}[\log G(z)]} \right)\\
=&\frac{\partial}{\partial t}\left(\arg \frac{\partial \log
G(re^{it})}{\partial t}\right).
\end{align*} Then
the assumptions imply that $\log G$ is convex if and only if $\log
F$ is convex.
\epf
\section{Goodman
and Saff's conjecture in $\mathcal{L}_p\mathcal{H}(G)$}\label{csw-sec5}

It is well-known that if an analytic function maps $\mathbb{D}$
univalently onto a convex domain, then it also maps each concentric
subdisk onto a convex domain (cf. \cite{pd}). Goodman and Saff
(\cite{a}) constructed an example of a function convex in the
vertical direction whose restriction to the disk $\mathbb{D}_{r}=\{z:\;|z| < r\}$
does not have that property for any radius $r$ in the interval
$\sqrt{2}-1 < r < 1$. In the same paper, they conjectured that the
radius $\sqrt{2}-1$ is best possible.

\bdefe\label{de1} A domain $D$ is {\it convex} in the direction
$e^{i\phi}$, if for every fixed complex number $z$, the set $D \cap
\{z + te^{i\phi} : t\in \mathbb{R}\}$ is either connected or empty. \edefe

Let $\mathcal{K}(\phi)$ (resp. $\mathcal{K}_{H}(\phi)$) denote the class of all
complex-valued analytic (resp. harmonic) univalent functions $f$ in
$\mathbb{D}$ with $f (\mathbb{D})$ convex in the direction
$e^{i\phi}$. If $f\in \mathcal{K}(\phi)$ (resp. $\mathcal{K}_{H}(\phi)$) is such that $f
(\mathbb{D})$ is convex in every direction (i.e. $f (\mathbb{D})$ is
a convex domain), then in this case we say that $f\in \mathcal{K}$ (resp.
$\mathcal{K}_{H}$). Ruscheweyh and Salinas \cite[Theorem $1$]{sl} ultimately
succeeded in proving the Goodman-Saff conjecture by showing that if
$f\in \mathcal{K}_{H}(\phi)$ and $r\in(0 ,\sqrt{2} - 1]$, then one has $f (rz)
\in \mathcal{K}_{H}(\phi)$. In particular, this gives the following.

\begin{Thm} $($\cite[Theorem $1$]{sl}$)$  Let $f\in \mathcal{K}_{H}$. Then for any  $r\in (0, \sqrt{2}-1]$, $f (rz) \in \mathcal{K}_{H}$.\end{Thm}

In view of the development of logharmonic
mappings, it is interesting to ask whether the same conjecture holds
for $\log$-polyharmonic mappings. Our result is as
follows.

\bthm\label{thm5} Suppose that $F\in \mathcal{L}_p\mathcal{H}(G)$ with $$F(z)
=f(z)h(\overline{z})\prod _{k=1}^{p}G(z)^{\lambda_{k}|z|^{2(k-1)}},$$
where $f$ and $h$ are constant functions,
 $\log F$ is univalent mapping, $\log G$ is univalent and convex, and $\mathcal{L}[\log G(z)]\neq0$, $\sum_{k=1}^{p}\lambda_{k}|z|^{2(k-1)}\neq0$ for $z\in \mathbb{D}\backslash \{0\}$.
Then $\log F$ sends the subdisk $\mathbb{D}_{r}$ onto a convex region for $r\in (0, \sqrt{2}-1]$.
\ethm
\bpf
Since $\log G$ is harmonic and if we further assume that it is
convex in $0<r\leq r_{0}=\sqrt{2}-1$, then by Theorem \ref{cor5}
and Theorem A, we have that $\log F$ is also convex in $0<r\leq
r_{0}=\sqrt{2}-1$. The proof is complete. \epf

\bigskip
\noindent {\bf Acknowledgements:}
The research was supported by NSF of
China (No. 11571216 and No. 11671127), NSF of Guangdong Province (No. 2014A030313471), and Project of ISTCIPU in Guangdong Province (No. 2014KGJHZ007).

The authors thank the referee very much for his/her careful reading of this paper and many useful suggestions.

\end{document}